\documentclass[12pt,a4paper]{amsart}
\usepackage{amsfonts}
\usepackage{amssymb}
\usepackage{amsmath}
\usepackage{amsthm}
\usepackage{cite}
\usepackage{graphicx}
\usepackage{amscd}
\usepackage{color}
\newtheorem{theorem}{Theorem}
\theoremstyle{plain}

\newtheorem{corollary}[theorem]{Corollary}
\newtheorem{lemma}[theorem]{Lemma}

\theoremstyle{definition}
\newtheorem{example}{Example}

\numberwithin{equation}{section}

\def\bi{\begin{itemize}}
\def\ei{\end{itemize}}

\newcommand{\R}{\mathbb{R}}
\newcommand{\C}{\mathbb{C}}

\renewcommand{\phi}{\varphi}

\DeclareMathOperator*{\esssup}{ess\,sup}

\newcommand{\supp}{\text{\rm supp}}

\newcommand{\bspm}{\left(\begin{smallmatrix}}\newcommand{\espm}{\end{smallmatrix}\right)}

\def\blem{\begin{lemma}}\def\elem{\end{lemma}}
\def\bthm{\begin{theorem}}\def\ethm{\end{theorem}}
\def\bcor{\begin{corollary}}\def\ecor{\end{corollary}}
\def\beq{\begin{equation}}\def\eeq{\end{equation}}

\begin{document}

\title[Interfaces for SGS Existence in 1D NLS]{Interfaces Supporting Surface Gap Soliton Ground States in the 1D Nonlinear Schr\"odinger Equation} 

\author{Tom\'{a}\v{s} Dohnal}
\address{T. Dohnal \hfill\break 
Department of Mathematics, Technical University Dortmund\hfill\break
D-44221 Dortmund, Germany}
\email{tomas.dohnal@math.tu-dortmund.de}

\author{Kaori Nagatou}
\address{K. Nagatou \hfill\break 
Institut f\"ur Analysis, Karlsruhe Institute of Technology (KIT) \hfill\break
D-76128 Karlsruhe, Germany\hfill\break 
and\hfill\break
Faculty of Science and Engineering, Waseda University\hfill \break 
3-4-1 Okubo, Shinjuku-ku, Tokyo 169-8555, Japan}
\email{kaori.nagatou@kit.edu}

\author{Michael Plum}
\address{M. Plum \hfill\break 
Institut f\"ur Analysis, Karlsruhe Institute of Technology (KIT)\hfill\break
D-76128 Karlsruhe, Germany}
\email{michael.plum@kit.edu}

\author{Wolfgang Reichel}
\address{W. Reichel \hfill\break 
Institut f\"ur Analysis, Karlsruhe Institute of Technology (KIT) \hfill\break
D-76128 Karlsruhe, Germany}
\email{wolfgang.reichel@kit.edu}

\date{\today}

\subjclass[2000]{Primary: 35Q55, 78M30, 65G20; Secondary: 35J20, 35Q60}
\keywords{nonlinear Schr\"odinger equation, surface gap soliton, ground state, variational methods, interface, periodic material, verified numerical enclosures}

\begin{abstract} 
We consider the problem of verifying the existence of $H^1$ ground states of the 1D nonlinear Schr\"odinger equation for an interface of two periodic structures:
$$-u'' +V(x)u -\lambda  u = \Gamma(x) |u|^{p-1}u \ \mbox{ on } \R$$
with $V(x) = V_1(x), \Gamma(x)=\Gamma_1(x)$ for $x\geq 0$ and $V(x) = V_2(x), \Gamma(x)=\Gamma_2(x)$ for $x<0$. Here $V_1,V_2,\Gamma_1,\Gamma_2$ are periodic, $\lambda <\min\sigma(-\tfrac{d^2}{dx^2}+V)$, and $p>1$. The article [T. Dohnal, M. Plum and W. Reichel, ``Surface Gap Soliton Ground States for the Nonlinear Schr\"odinger Equation,'' \textit{Comm. Math. Phys.} \textbf{308}, 511-542 (2011)] provides in the 1D case an existence criterion in the form of an integral inequality involving the linear potentials $V_{1},V_2$ and the Bloch waves of the operators $-\tfrac{d^2}{dx^2}+V_{1,2}-\lambda$. We choose here the classes of piecewise constant and piecewise linear potentials $V_{1,2}$ and check this criterion for a set of parameter values. In the piecewise constant case  the Bloch waves are calculated explicitly and in the 
piecewise linear case verified enclosures of the Bloch waves are computed numerically. The integrals in the criterion are evaluated via interval arithmetic so that rigorous existence statements are produced. Examples of interfaces supporting ground states are reported including such, for which ground state existence follows for all periodic $\Gamma_ {1,2}$ with $\esssup \Gamma_{1,2}>0$.
\end{abstract}

\maketitle


\section{Introduction}

An interface between two nonlinear periodic media in the $n-$dimensional nonlinear Schr\"odinger model can act as a waveguide so that localized solutions, so called surface gap solitons (SGS), exist as shown analytically in \cite{DPR11}. Experimentally such waveguiding has been demonstrated in nonlinear photonic crystals, see e.g. \cite{Rosberg06,Suntsov08,Szameit07}. There are also a number of numerical observations of SGSs in the 1D and 2D nonlinear Schr\"odinger equation (NLS), see e.g. \cite{BD11,DP08,Kartashov06,MHChSS06}. 

In \cite{DPR11} an existence criterion for strong ground states of the $n-$dimensional NLS 
$$
(-\Delta +V(x)-\lambda)u=\Gamma(x) |u|^{p-1}u, \ x\in \R^n  \leqno{\text{($n$-NLS)}}
$$
was proved with $V(x) = V_1(x), \Gamma(x)=\Gamma_1(x)$ for $x_1\geq 0$ and $V(x) = V_2(x), \Gamma(x)=\Gamma_2(x)$ for $x_1<0$ under the condition $\lambda < \min \sigma(-\Delta + V)$. The functions $V_1, V_2, \Gamma_1, \Gamma_2$ are assumed periodic in each coordinate direction and the exponent $p$ satisfies $p\in(1,2^\ast)$, where $2^\ast=\frac{2n}{n-2}$ if $n\geq 3$ and $2^\ast=\infty$ if $n=1,2$. A strong ground state is defined to be a minimizer of the corresponding total energy $\int_{\R^n} \tfrac{1}{2}(|\nabla u|^2+(V(x)-\lambda)u^2)-\tfrac{1}{p+1}\Gamma(x)|u|^{p+1}dx$ restricted to the Nehari manifold $N=\{u\in H^1(\R^n)\setminus\{0\}: \int_{\R^n}|\nabla u|^2+(V(x)-\lambda)u^2 -\Gamma(x)|u|^{p+1}dx=0\}$. The results of \cite{DPR11} include sufficient conditions for the existence of strong ground states. These conditions involve information about the strong ground states\footnote{The existence of strong ground states of the purely periodic problem on $\R^n$ was proved in \cite{Pankov05}.} $w_1, w_2$ of the purely periodic problems ($n$-NLS) with 
$V=V_1,\Gamma=\Gamma_1$ on $\R^n$ and $V=V_2,\Gamma=\Gamma_2$ on $\R^n$ respectively. In the case $n=1$ 
these conditions could be formulated in terms of the Bloch waves of the two purely periodic linear problems. As neither the ground states $w_1, w_2$ nor the Bloch waves are generally known explicitly, \cite{DPR11} did not produce explicit examples of ground state supporting interfaces except for an example where the potentials are related by scaling: $V_1(x) = k^2 V_2(k x)$, $\Gamma_1(x) = \gamma^2 \Gamma_2(k x)$ with certain conditions on $k$ and $\gamma$, see Theorem 5 in  \cite{DPR11}. All other existence examples were asymptotic; either in $\lambda$ or in $\Gamma_1-\Gamma_2$.

The most practical existence criteria in  \cite{DPR11} are those for the 1D case $n=1$. In this article we provide a number of explicit 1D examples of interfaces satisfying these criteria. We consider, therefore
\begin{equation}
-u'' +V(x)u -\lambda  u = \Gamma(x) |u|^{p-1}u \mbox{ on } \R
\label{int_nls}
\end{equation}
with
\beq
V(x) = \left\{
\begin{array}{ll}
V_1(x), & x\geq 0, \vspace{\jot}\\
V_2(x), & x<0,
\end{array} \right.
\label{E:V_interf}
\eeq
and
\beq
\Gamma(x) = \left\{
\begin{array}{ll}
\Gamma_1(x), & x\geq 0, \vspace{\jot}\\
\Gamma_2(x), & x<0
\end{array} \right.
\label{E:Gamma_interf}
\eeq
under the assumptions
\begin{itemize}
\item[(H1)] $V_1, V_2, \Gamma_1, \Gamma_2$ are $1$-periodic,
\item[(H2)] $\esssup \Gamma_i>0$, $i=1,2$,
\item[(H3)] $1<p<\infty$,
\item[(H4)] $\lambda< \min \sigma(-\frac{d^2}{dx^2} +V)$,
\end{itemize}
which were needed in  \cite{DPR11}. 

Next, recall the criterion given in Theorem 7 in \cite{DPR11} for the existence of SGS ground states of \eqref{int_nls}.
\begin{theorem}
Assume (H1)--(H4) and for $i=1,2$ define by $c_i$ the energy  of a strong ground state of $(-\frac{d^2}{dx^2} +V_i(x) -\lambda)u=\Gamma_{i}(x)|u|^{p-1}u$ on $\R$. 
\begin{itemize}
\item[(a)] If $c_1\leq c_2$, then a sufficient condition for the existence of a strong ground state of \eqref{int_nls} is
\begin{equation}
I_1:=\int_{-1}^0 \Big(V_2(x)-V_1(x)\Big) u^{(1)}_-(x)^2\,dx<0,
\label{gen_cond1}
\end{equation}
where $u^{(1)}_-(x)=p^{(1)}_-(x)e^{\kappa_1 x}$, with $\kappa_1>0$ and $p_-^{(1)}$ 1-periodic, is the Bloch mode decaying at $-\infty$ of $-\frac{d^2}{dx^2}+ V_1(x)-\lambda$.
\item[(b)] If $c_1\geq c_2$, then a sufficient condition for the existence of a strong ground state of \eqref{int_nls} is
\begin{equation}
I_2:=\int_{0}^1 \Big(V_1(x)-V_2(x)\Big) u^{(2)}_+(x)^2\,dx<0,
\label{gen_cond2}
\end{equation}
where $u^{(2)}_+(x)=p^{(2)}_+(x)e^{-\kappa_2 x}$, with $\kappa_2>0$ and $p_+^{(2)}$ 1-periodic, is the Bloch mode decaying at $+\infty$ of $-\frac{d^2}{dx^2}+ V_2(x)-\lambda$.
\label{T:Bloch_cond}
\end{itemize}
\end{theorem}
When the ordering of $c_1, c_2$ is unknown, Theorem \ref{T:Bloch_cond} can still be used by establishing negativity of \textit{both} of the integrals $I_1$ and $I_2$:
\bcor\label{C:both_neg}
If both $I_1,I_2<0$, then a strong ground state exists irrespectively of the order of $c_1$ and $c_2$, and thus, of the choice of $\Gamma_1,\Gamma_2$ and $p$ (within the assumptions (H1)-(H3)).
\ecor
As seen from \eqref{gen_cond1} and \eqref{gen_cond2}, any ordering $V_1(x) \geq V_2(x)$ or $V_1(x) \leq V_2(x)$ on $[0,1]$ leads to $I_1\leq 0\leq I_2$ or $I_2\leq 0\leq I_1$, whence the assumption of corollary \ref{C:both_neg} is not satisfied.    

If information on the ordering of the ground state energies $c_1,c_2$ is available, the corresponding criterion (a) or (b) in Theorem \ref{T:Bloch_cond} can be checked. This is the case, for instance, with the dislocation interface
\beq\label{E:disloc_interf}
V(x) = \left\{
\begin{array}{ll}
V_0(x+\tau_1), & x\geq 0, \vspace{\jot}\\
V_0(x-\tau_2), & x<0,
\end{array} \right.
\qquad
\Gamma(x) = \left\{
\begin{array}{ll}
\Gamma_0(x+\tau_1), & x\geq 0, \vspace{\jot}\\
\Gamma_0(x-\tau_2), & x<0,
\end{array} \right.
\eeq
where $V_0,\Gamma_0$ are $1-$periodic, $\esssup \Gamma_0>0$, and where $\tau_1,\tau_2\in [0,1]$ are the dislocation parameters. In this case $c_1=c_2$ so that we have
\bcor
For the interface \eqref{E:disloc_interf} suppose $I_1<0$ or $I_2<0$. Then there exists a strong ground state irrespectively of the choice of $p$ (within (H3)) and of $\Gamma_0$ (within the above assumptions).
\ecor

We will use direct constructive approaches to verify the respective
conditions \eqref{gen_cond1} and \eqref{gen_cond2}. These require mainly the Bloch waves of the two purely periodic linear problems on $\R$. We consider two types of potentials $V_{1,2}$: piecewise constant and piecewise linear. For piecewise constant potentials we
calculate the needed Bloch modes in closed form by hand, so that we can check these
conditions directly. For piecewise linear potentials we use a computer-assisted approach,
i.e. we compute verified enclosures for the Bloch modes, and use these to enclose the
integrals $I_1$ and $I_2$ in \eqref{gen_cond1} and \eqref{gen_cond2}. These computer-assisted results are completely verified and thus give a rigorous mathematical proof since all numerical errors are
taken into account. In principle, even much more general potentials can be treated by this approach; we have chosen piecewise linear ones for simplicity.

In the rest of the paper we impose the condition 
\beq\label{E:lam_cond}
\lambda < \min \{\inf{V_1},\inf{V_2}\},
\eeq
which ensures (H4) without having to actually calculate the spectrum.

In Section \ref{S:pwise_const} below we consider interfaces with piecewise constant potentials  $V_1$ and $V_2$  and in Section \ref{S:pwise_lin} we study the case of piecewise linear potentials.

\section{Piecewise Constant Potentials $V_1$ and $V_2$}\label{S:pwise_const}

When the potentials $V_1$ and $V_2$ are piecewise constant, the integrals $I_1, I_2$ can be calculated explicitly although the closed form involves the inverse of a transcendental function. We calculate the formulas for $I_1, I_2$ explicitly and evaluate these numerically for a set of parameter values. The evaluation is done in interval arithmetic (using the Matlab toolbox Intlab \cite{Ru99a}). The resulting values of $I_1, I_2$ are thus enclosed in intervals and when the supremum of such an interval is negative, the corresponding integral $I_1$ or $I_2$ is then verified to be negative.

\subsection{Bloch Waves for a Piecewise Constant Potential}

\bi
\item[]
\ei
\vspace{-0.2cm}
Let
\beq\label{E:pw_const_V}
V_0(x) = \left\{
\begin{array}{ll}
a, & 0\leq x<s, \vspace{\jot}\\
b, & s\leq x <1
\end{array} \right.
\eeq
with $a, b \in \R, a \neq b$ and $s \in (0,1)$.

The Bloch waves of $-u'' +V_0(x)u=\lambda u$ on $\R$ have the form
$$u_\pm(x) = p_{\pm}(x) e^{\mp \kappa x} \quad \text{with} \ p_{\pm}(x+1) = p_{\pm}(x), \kappa > 0$$
since $\lambda$ lies in the resolvent set of $-\tfrac{d^2}{dx^2}+V_0$, see \cite{Eastham}.
At the same time, due to the piecewise constant nature of $V_0$
\beq\label{E:Bloch_wvs}
u_\pm(x) = \left\{
\begin{array}{ll}
\xi_1^{\pm}e^{x\sqrt{a-\lambda}} + \xi_2^{\pm}e^{-x\sqrt{a-\lambda}}, & 0 \leq x < s, \vspace{\jot}\\
\xi_3^{\pm}e^{x\sqrt{b-\lambda}} + \xi_4^{\pm}e^{-x\sqrt{b-\lambda}}, & s \leq x<1.
\end{array} \right.
\eeq
Note that $\lambda < \min\{a,b\}$ due to \eqref{E:lam_cond}. The vectors $\xi^{\pm}$ are determined via the $C^1$ condition for $u_\pm(x)$ at $x=s$ and the condition that the Floquet multipliers of $u_\pm$ are $e^{\mp \kappa}$ respectively. For $\xi^+$ we thus obtain the system
$$
A(\kappa)\xi^+ = 0, \qquad A(\kappa) = \bspm
e^{s\sqrt{a-\lambda}} & e^{-s\sqrt{a-\lambda}} & -e^{s\sqrt{b-\lambda}} & -e^{-s\sqrt{b-\lambda}}\\
\sqrt{a-\lambda}~e^{s\sqrt{a-\lambda}} & -\sqrt{a-\lambda}~e^{-s\sqrt{a-\lambda}} & -\sqrt{b-\lambda}~e^{s\sqrt{b-\lambda}} & \sqrt{b-\lambda}~e^{-s\sqrt{b-\lambda}}\\
-e^{-\kappa} &  -e^{-\kappa} & e^{\sqrt{b-\lambda}}& e^{-\sqrt{b-\lambda}}\\
-\sqrt{a-\lambda}~e^{-\kappa} &  \sqrt{a-\lambda}~e^{-\kappa} & \sqrt{b-\lambda}~e^{\sqrt{b-\lambda}} & -\sqrt{b-\lambda}~e^{-\sqrt{b-\lambda}}\\
\espm.
$$
Solving $\det(A(\kappa))=0$ yields
\beq\label{E:kappa}
\begin{split}
\cosh(\kappa) = \frac{1}{4\sqrt{a-\lambda}\sqrt{b-\lambda}}&\left[(\sqrt{a-\lambda}+\sqrt{b-\lambda})^2\cosh\left(s\sqrt{a-\lambda} + (1-s)\sqrt{b-\lambda}\right)\right.\\
&\left. -(\sqrt{a-\lambda}-\sqrt{b-\lambda})^2\cosh\left(s\sqrt{a-\lambda} - (1-s)\sqrt{b-\lambda}\right)\right].
\end{split}
\eeq
The solution vector $\xi^+$ is proportional to
\beq
\label{E:xi_plus}
\begin{split}
\xi_1^+ & = \frac{e^{-s\sqrt{a-\lambda}}}{2\left(e^{s(\sqrt{a-\lambda}+\sqrt{b-\lambda})-\kappa}-e^{\sqrt{b-\lambda}}\right)}\left[4e^{s\sqrt{a-\lambda}-\kappa}\sqrt{a-\lambda}\sqrt{b-\lambda} \right. \\
& \hspace{0.5cm}\left. - \left(\sqrt{a-\lambda}+\sqrt{b-\lambda}\right)^2e^{(s-1)\sqrt{b-\lambda}}
+\left(\sqrt{a-\lambda}-\sqrt{b-\lambda}\right)^2 e^{(1-s)\sqrt{b-\lambda}} \right],\\
\xi_2^+ &= (a-b)\frac{\sinh((1-s)\sqrt{b-\lambda})}{e^{s\sqrt{b-\lambda}-\kappa} - e^{\sqrt{b-\lambda}-s\sqrt{a-\lambda}}},\\
\xi_3^+ &= \left(a-\lambda+\sqrt{a-\lambda}\sqrt{b-\lambda}\right)\frac{e^{(\sqrt{a-\lambda}-\sqrt{b-\lambda})s-\kappa}-e^{-\sqrt{b-\lambda}}}{e^{(\sqrt{a-\lambda}+\sqrt{b-\lambda})s-\kappa}-e^{\sqrt{b-\lambda}}},\\
\xi_4^+ &= \lambda -a + \sqrt{a-\lambda}\sqrt{b-\lambda}.
\end{split}
\eeq
The system for $\xi^-$ reads $A(-\kappa)\xi^-=0$, so that
\beq
\label{E:xi_min}
\xi^-(\kappa) = \xi^+(-\kappa).
\eeq

\subsection{The Dislocation Interface}

\bi
\item[]
\ei
\vspace{-0.2cm}
Let us consider the dislocation interface \eqref{E:disloc_interf} with $\tau_1=\tau_2=:\tau$ for the piecewise constant potential
$$
V_0(x) = \left\{
\begin{array}{ll}
a, & 0\leq x<1/2, \vspace{\jot}\\
b, & 1/2 \leq x <1
\end{array} \right. 
$$
with $a,b\in \R, a\neq b$. In this case $u_\pm^{(1)}(x)=u_\pm(x+\tau)$ and $u_\pm^{(2)}(x)=u_\pm(x-\tau)$ with $u_\pm$ in \eqref{E:Bloch_wvs}, $\xi^{\pm}$ in \eqref{E:xi_plus}, \eqref{E:xi_min}  and $\kappa$ in \eqref{E:kappa}, where we set $s=1/2$.

A direct integration then produces for $0<\tau<1/4$
\beq
\begin{split}
I_1(\tau)&= \frac{b-a}{2\sqrt{a-\lambda}}\left[\xi_1^{-^2}\left(e^{-2\kappa+4\tau\sqrt{a-\lambda}}-e^{2\tau\sqrt{a-\lambda}}(e^{-2 \kappa}-1)-1\right) \right.\\
                                      & \hspace{2cm}\left.- \xi_2^{-^2}\left(e^{-2\kappa-4\tau\sqrt{a-\lambda}}-e^{-2\tau\sqrt{a-\lambda}}(e^{-2 \kappa}-1)-1\right) \right]\\ 
& + \frac{a-b}{2\sqrt{b-\lambda}}e^{-2 \kappa}\left[\xi_3^{-^2}\left(e^{(1+4\tau)\sqrt{b-\lambda}}-e^{\sqrt{b-\lambda}}\right)-\xi_4^{-^2}\left(e^{-(1+4\tau)\sqrt{b-\lambda}}-e^{-\sqrt{b-\lambda}}\right)\right] \\
& + 2\tau(b-a)\left[(e^{-2\kappa}+1)\xi_1^-\xi_2^- -2e^{-2\kappa}\xi_3^-\xi_4^-\right] ,
\end{split}
\eeq
for $1/4\leq \tau<1/2$
\beq
\begin{split}
I_1(\tau)&= \frac{b-a}{2\sqrt{a-\lambda}}\left[\xi_1^{-^2}\left(e^{-2\kappa+\sqrt{a-\lambda}}-e^{2\tau\sqrt{a-\lambda}}(e^{-2 \kappa}-1)-e^{(4\tau-1)\sqrt{a-\lambda}}\right) \right.\\
                                      & \hspace{2cm}\left.- \xi_2^{-^2}\left(e^{-2\kappa-\sqrt{a-\lambda}}-e^{-2\tau\sqrt{a-\lambda}}(e^{-2 \kappa}-1)-e^{-(4\tau-1)\sqrt{a-\lambda}}\right)\right]\\ 
& + \frac{a-b}{2\sqrt{b-\lambda}}e^{-2 \kappa}\left[\xi_3^{-^2}\left(e^{2\sqrt{b-\lambda}}-e^{4\tau\sqrt{b-\lambda}}\right)-\xi_4^{-^2}\left(e^{-2\sqrt{b-\lambda}}-e^{-4\tau\sqrt{b-\lambda}}\right)\right] \\
& + (1-2\tau)(b-a)\left[(e^{-2\kappa}+1)\xi_1^-\xi_2^- - 2e^{-2\kappa}\xi_3^-\xi_4^-\right] ,
\end{split}
\eeq
for $1/2\leq \tau<3/4$
\beq
\begin{split}
I_1(\tau)&= \frac{b-a}{2\sqrt{a-\lambda}}\left[\xi_1^{-^2}\left(e^{(4\tau-2)\sqrt{a-\lambda}}-1\right)-\xi_2^{-^2}\left(e^{-(4\tau-2)\sqrt{a-\lambda}}-1\right)\right]\\
& + \frac{a-b}{2\sqrt{b-\lambda}}\left[\xi_3^{-^2}\left(e^{-2\kappa+(4\tau-1)\sqrt{b-\lambda}}-e^{2\tau\sqrt{b-\lambda}}(e^{-2 \kappa}-1)-e^{\sqrt{b-\lambda}}\right) \right.\\
                                      & \hspace{2cm}\left.- \xi_4^{-^2}\left(e^{-2\kappa-(4\tau-1)\sqrt{b-\lambda}}-e^{-2\tau\sqrt{b-\lambda}}(e^{-2 \kappa}-1)-e^{-\sqrt{b-\lambda}}\right)\right]\\ 
& + (2\tau-1)(a-b)\left[(e^{-2\kappa}+1)\xi_3^-\xi_4^- - 2\xi_1^-\xi_2^-\right] ,
\end{split}
\eeq
and for $3/4\leq \tau<1$
\beq
\begin{split}
I_1(\tau)&= \frac{b-a}{2\sqrt{a-\lambda}}\left[\xi_1^{-^2}\left(e^{\sqrt{a-\lambda}}-e^{(4\tau-3)\sqrt{a-\lambda}}\right)-\xi_2^{-^2}\left(e^{-\sqrt{a-\lambda}}-e^{-(4\tau-3)\sqrt{a-\lambda}}\right)\right]\\
& + \frac{a-b}{2\sqrt{b-\lambda}}\left[\xi_3^{-^2}\left(e^{-2\kappa+2\sqrt{b-\lambda}}-e^{2\tau\sqrt{b-\lambda}}(e^{-2 \kappa}-1)-e^{(4\tau-2)\sqrt{b-\lambda}}\right) \right.\\
                                      & \hspace{2cm}\left.- \xi_4^{-^2}\left(e^{-2\kappa-2\sqrt{b-\lambda}}-e^{-2\tau\sqrt{b-\lambda}}(e^{-2 \kappa}-1)-e^{-(4\tau-2)\sqrt{b-\lambda}}\right)\right]\\ 
& + (2-2\tau)(a-b)\left[(e^{-2\kappa}+1)\xi_3^-\xi_4^- - 2\xi_1^-\xi_2^-\right].
\end{split}
\eeq

For the integral $I_2$ we have for $0<\tau<1/4$
\beq
\begin{split}
I_2(\tau)&= \frac{b-a}{2\sqrt{a-\lambda}}\left[\xi_1^{+^2}\left(e^{\sqrt{a-\lambda}}-e^{(1-4\tau)\sqrt{a-\lambda}}\right)-\xi_2^{+^2}\left(e^{-\sqrt{a-\lambda}}-e^{-(1-4\tau)\sqrt{a-\lambda}}\right)\right]\\
& + \frac{a-b}{2\sqrt{b-\lambda}}\left[\xi_3^{+^2}\left(e^{2(1-\tau)\sqrt{b-\lambda}}(1-e^{2 \kappa})-e^{(2-4\tau)\sqrt{b-\lambda}}+e^{2\sqrt{b-\lambda}+2\kappa}\right) \right.\\
                                      & \hspace{2cm}\left.- \xi_4^{+^2}\left(e^{-2(1-\tau)\sqrt{b-\lambda}}(1-e^{2\kappa})-e^{-(2-4\tau)\sqrt{b-\lambda}}+e^{-2\sqrt{b-\lambda}+2\kappa}\right)\right]\\ 
& + 2\tau(a-b)\left[(e^{2\kappa}+1)\xi_3^+\xi_4^+ - 2\xi_1^+\xi_2^+\right],
\end{split}
\eeq
for $1/4\leq\tau<1/2$
\beq
\begin{split}
I_2(\tau)&= \frac{b-a}{2\sqrt{a-\lambda}}\left[\xi_1^{+^2}\left(e^{(2-4\tau)\sqrt{a-\lambda}}-1\right)-\xi_2^{+^2}\left(e^{-(2-4\tau)\sqrt{a-\lambda}}-1\right)\right]\\
& + \frac{a-b}{2\sqrt{b-\lambda}}\left[\xi_3^{+^2}\left(e^{2(1-\tau)\sqrt{b-\lambda}}(1-e^{2 \kappa})-e^{\sqrt{b-\lambda}}+e^{(3-4\tau)\sqrt{b-\lambda}+2\kappa}\right) \right.\\
                                      & \hspace{2cm}\left.- \xi_4^{+^2}\left(e^{-2(1-\tau)\sqrt{b-\lambda}}(1-e^{2 \kappa})-e^{-\sqrt{b-\lambda}}+e^{-(3-4\tau)\sqrt{b-\lambda}+2\kappa}\right)\right]\\ 
& + (1-2\tau)(a-b)\left[(e^{2\kappa}+1)\xi_3^+\xi_4^+ - 2\xi_1^+\xi_2^+\right],
\end{split}
\eeq
for $1/2\leq\tau<3/4$
\beq
\begin{split}
I_2(\tau)&=  \frac{b-a}{2\sqrt{a-\lambda}}\left[\xi_1^{+^2}\left(e^{(2-2\tau)\sqrt{a-\lambda}}(1-e^{2 \kappa})-e^{(3-4\tau)\sqrt{a-\lambda}}+e^{\sqrt{a-\lambda}+2\kappa}\right) \right.\\
                                      & \hspace{2cm}\left.- \xi_2^{+^2}\left(e^{-(2-2\tau)\sqrt{a-\lambda}}(1-e^{2 \kappa})-e^{-(3-4\tau)\sqrt{a-\lambda}}+e^{-\sqrt{a-\lambda}+2\kappa}\right)\right]\\ 
& + \frac{a-b}{2\sqrt{b-\lambda}}e^{2 \kappa}\left[\xi_3^{+^2}\left(e^{2\sqrt{b-\lambda}}-e^{(4-4\tau)\sqrt{b-\lambda}}\right)-\xi_4^{+^2}\left(e^{-2\sqrt{b-\lambda}}-e^{-(4-4\tau)\sqrt{b-\lambda}}\right)\right] \\
& + (2\tau-1)(b-a)\left[(e^{2\kappa}+1)\xi_1^+\xi_2^+ - 2e^{2\kappa}\xi_3^+\xi_4^+\right] ,
\end{split}
\eeq
and for $3/4\leq\tau<1$
\beq
\begin{split}
I_2(\tau)&=  \frac{b-a}{2\sqrt{a-\lambda}}\left[\xi_1^{+^2}\left(e^{(2-2\tau)\sqrt{a-\lambda}}(1-e^{2 \kappa})-1+e^{(4-4\tau)\sqrt{a-\lambda}+2\kappa}\right) \right.\\
                                      & \hspace{2cm}\left.- \xi_2^{+^2}\left(e^{-(2-2\tau)\sqrt{a-\lambda}}(1-e^{2 \kappa})-1+e^{-(4-4\tau)\sqrt{a-\lambda}+2\kappa}\right)\right]\\ 
& + \frac{a-b}{2\sqrt{b-\lambda}}e^{2 \kappa}\left[\xi_3^{+^2}\left(e^{(5-4\tau)\sqrt{b-\lambda}}-e^{\sqrt{b-\lambda}}\right)-\xi_4^{+^2}\left(e^{-(5-4\tau)\sqrt{b-\lambda}}-e^{-\sqrt{b-\lambda}}\right)\right] \\
& + (2-2\tau)(b-a)\left[(e^{2\kappa}+1)\xi_1^+\xi_2^+ - 2e^{2\kappa}\xi_3^+\xi_4^+\right].
\end{split}
\eeq

\begin{figure}[!ht]
  \begin{center}
  \scalebox{0.6}{\includegraphics{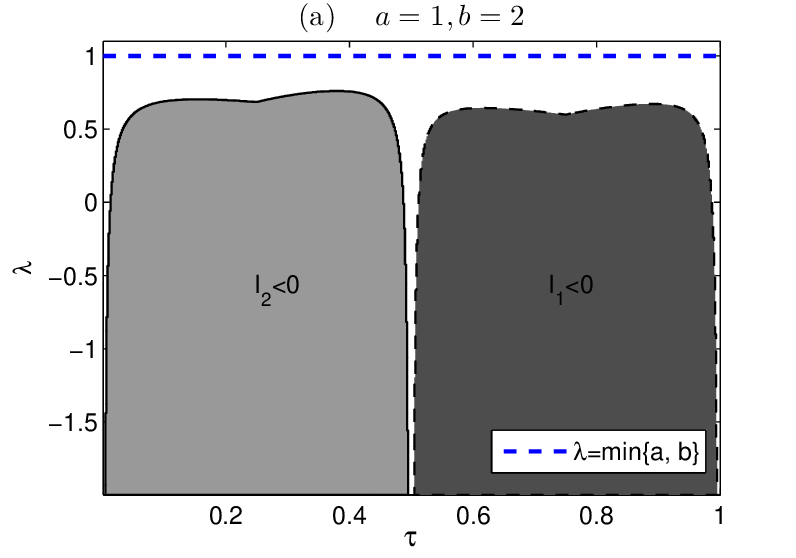}}
  \scalebox{0.6}{\includegraphics{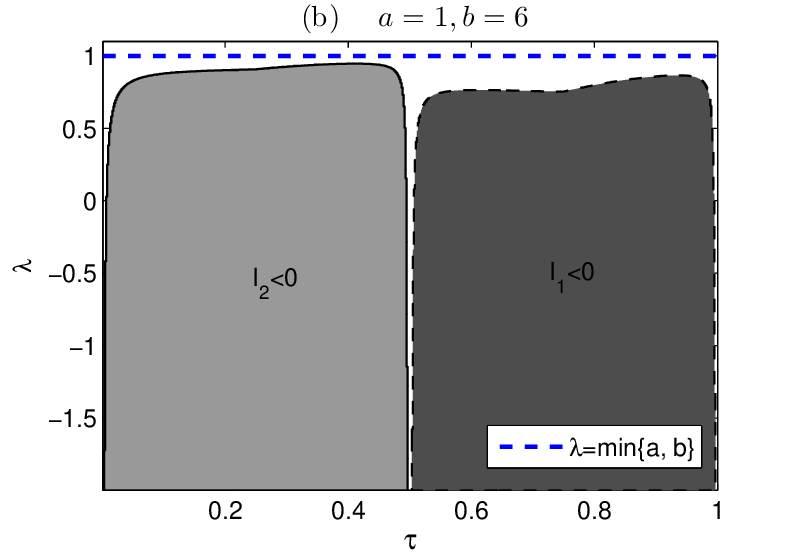}}
  \end{center}
  \caption{Regions in the $(\tau,\lambda)$ plane, where suprema of the enclosures of $I_1$ and $I_2$ are negative for the dislocation interface with a piecewise constant $V_0$ and $a=1$. (a) $b=2$; (b) $b=6$ (computations performed in interval arithmetic).}
  \label{F:disloc_2D}
\end{figure}

For the dislocation interfaces with $a=1$ and $b\in \{2,6\}$ Fig. \ref{F:disloc_2D} shows regions of the $(\tau,\lambda)$ plane where the integral $I_1$ or $I_2$ is negative, i.e. where ground state existence is guaranteed. These regions were computed using interval arithmetic. The domain $[0,1]\times [-2,0.98]$ in the $(\tau,\lambda)$-plane was completely covered\footnote{All intervals whose endpoints are not floating-point numbers are safely enclosed in slightly larger intervals by use of Intlab.}  by two dimensional intervals (squares) of size $1/600$ along each dimension and when for a given square the interval arithmetic evaluation produced $I_k<0$ (which means a negative supremum of the enclosure of $I_k$), the square was shaded.
Note that it is a priori clear that both $I_1$ and $I_2$ are zero at $\tau=0,1/2$, and $1$ because for these values $V_0(x+\tau)-V_0(x-\tau)\equiv 0$ due to the $1-$periodicity of $V_0$ and the integrands in $I_1, I_2$ thus vanish . The use of interval arithmetic in the computations then necessarily results in small neighborhoods of $\tau=0, 1/2, 1$ where the sign of the integrals cannot be determined. The quantities $I_1,I_2$ contain $\sqrt{a-\lambda}$ and $\sqrt{b-\lambda}$ in the denominator. To reduce the amount of round-off error for $\lambda$ close to $\min\{a,b\}$, we compute $\sqrt{a-\lambda}\sqrt{b-\lambda} \ I_{1,2}$ instead of $I_{1,2}$.

\begin{figure}[!ht]
  \begin{center}
  \scalebox{0.5}{\includegraphics{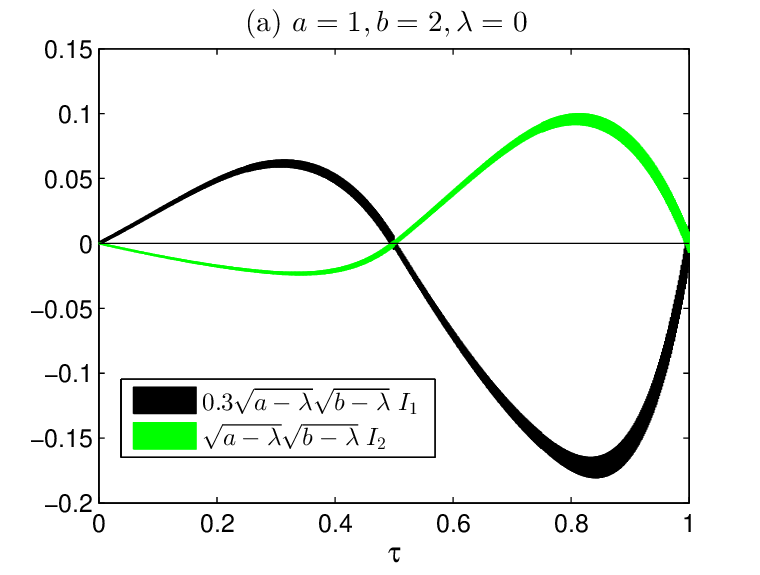}}
  \scalebox{0.5}{\includegraphics{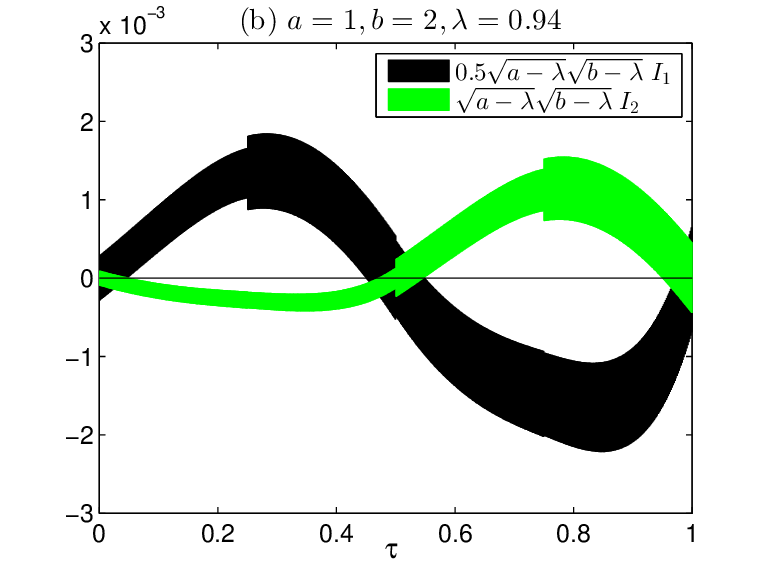}}
  \end{center}
  \caption{Scaled plots of the interval enclosures of $\sqrt{a-\lambda}\sqrt{b-\lambda} \ I_1$ and $\sqrt{a-\lambda}\sqrt{b-\lambda} \ I_2$ as functions of $\tau$ corresponding to Fig. \ref{F:disloc_2D} (a). In (a) $\lambda=0$ and in (b) $\lambda = 0.94$ (computations performed in interval arithmetic with the $\tau$-interval width $1/600$).}
  \label{F:disloc_tau_b2}
\end{figure}

\begin{figure}[!ht]
  \begin{center}
  \scalebox{0.55}{\includegraphics{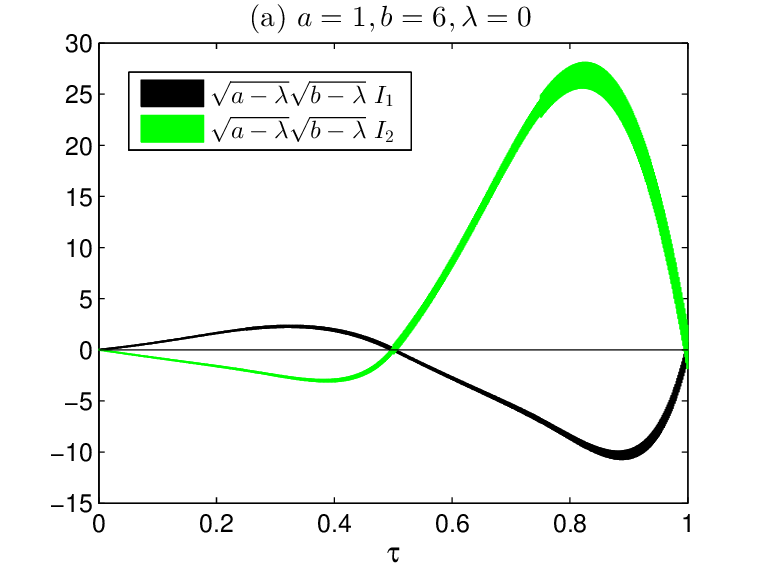}}
  \scalebox{0.55}{\includegraphics{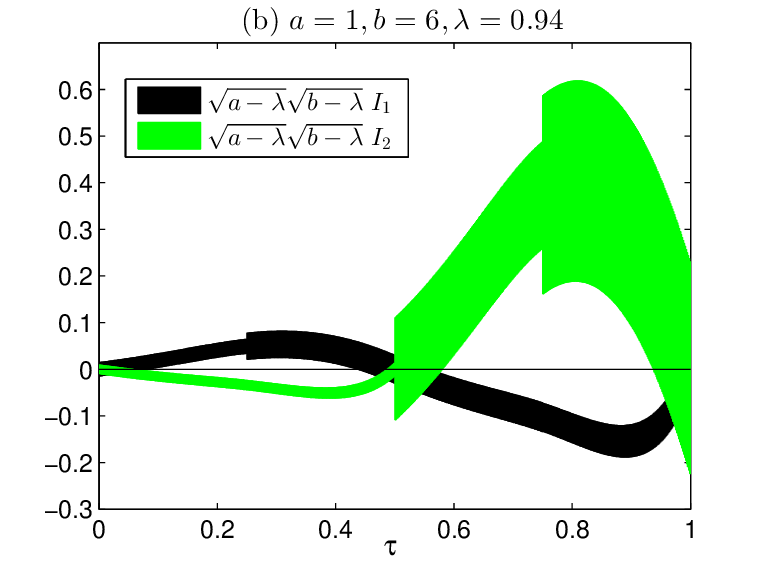}}
  \end{center}
  \caption{Scaled plots of the interval enclosures of $\sqrt{a-\lambda}\sqrt{b-\lambda} \ I_1$ and $\sqrt{a-\lambda}\sqrt{b-\lambda} \ I_2$ as functions of $\tau$ with $a$ and $b$ as in Fig. \ref{F:disloc_2D} (b). In (a) $\lambda=0$ and in (b) $\lambda = 0.94$ (computations performed in interval arithmetic with the $\tau$-interval width $1/600$).}
  \label{F:disloc_tau_b6}
\end{figure}

As Fig. \ref{F:disloc_2D} shows, ground state existence is guaranteed in the cases $b=2$ in almost the entire subset $(\tau,\lambda)\in [0,1]\times [-2,0.5]$ and in the case $b=6$  in almost the entire subset $(\tau,\lambda)\in [0,1]\times [-2,0.7]$ of the parameter domain. As the results in Fig. \ref{F:disloc_tau_b2} and \ref{F:disloc_tau_b6} show, these subsets can, in fact, be enlarged by reducing the interval size in interval arithmetic.

In Fig. \ref{F:disloc_tau_b2}  and \ref{F:disloc_tau_b6} the enclosures of $\sqrt{a-\lambda}\sqrt{b-\lambda} \ I_1$ and $\sqrt{a-\lambda}\sqrt{b-\lambda} \ I_2$ as functions of $\tau$ are plotted for the two examples in Fig. \ref{F:disloc_2D} at two values of $\lambda$, namely $\lambda =0$ and $\lambda =0.94$. As Fig.   \ref{F:disloc_tau_b2} (b) and \ref{F:disloc_tau_b6} (b) show, even at $\lambda =0.94$ there are $\tau$-intervals for which $I_1<0$ or $I_2<0$. This appears to contradict Fig. \ref{F:disloc_2D} but it is due to the use of two-dimensional intervals (intervals in the $(\tau,\lambda)$-plane) for input in Fig. \ref{F:disloc_2D}  and one-dimensional intervals in Fig. \ref{F:disloc_tau_b2}  and \ref{F:disloc_tau_b6}. The larger overestimation due to interval arithmetic in the case of two dimensional intervals leads to a smaller region where negativity of $I_1$ or $I_2$ is verified. 


\subsection{An interface where the is no knowledge about $c_1, c_2$}\label{S:pwise_const_gen_int}

\bi
\item[]
\ei
\vspace{-0.2cm}
Here we consider a general interface \eqref{E:V_interf}, \eqref{E:Gamma_interf} with the piecewise constant structure \eqref{E:pw_const_V} for both $V_1$ and $V_2$. In detail, $V_1$ is given by \eqref{E:pw_const_V} with $(a,b,s)$ replaced by $(a_1,b_1,s_1)$ and $V_2$ is given by \eqref{E:pw_const_V} with $(a,b,s)$ replaced by $(a_2,b_2,s_2)$. For simplicity we choose the jump locations in the middle of the periodicity cell: $s_1=s_2=1/2$. 

The Bloch waves $u^{(1)}_\pm$ are now given by \eqref{E:Bloch_wvs} and \eqref{E:kappa}-\eqref{E:xi_min} with $(a,b,s)$ replaced by $(a_1,b_1,1/2)$. We denote the resulting $\kappa$ in \eqref{E:kappa} by $\kappa_1$. Analogously we obtain  $u^{(2)}_\pm$ and denote the resulting vectors in \eqref{E:xi_plus},\eqref{E:xi_min} by $\zeta^\pm$.

Because the ordering of $c_1$ and $c_2$ is unknown in this case, Theorem \ref{T:Bloch_cond} can be used to prove ground state existence only if both $I_1$ and $I_2$ are negative. If this occurs, the existence of a strong ground state is then completely independent of the nonlinear periodic coefficients $\Gamma_1, \Gamma_2$ and of $p$ (within (H1)-(H3)). We show below that such cases occur. 

The integrals from Theorem \ref{T:Bloch_cond} now become
\beq
\begin{split}
I_1&= e^{-2\kappa_1}(a_2-a_1)\xi_1^-\xi_2^- + \frac{e^{-2\kappa_1}(a_2-a_1)}{2\sqrt{a_1-\lambda}}\left[\xi_1^{-^2}(e^{\sqrt{a_1-\lambda}}-1)+\xi_2^{-^2}(1-e^{-\sqrt{a_1-\lambda}})\right]\\
& +  e^{-2\kappa_1}(b_2-b_1)\xi_3^-\xi_4^- + \frac{e^{-2\kappa_1}(b_2-b_1)}{2\sqrt{b_1-\lambda}}\left[\xi_3^{-^2}(e^{2\sqrt{b_1-\lambda}}-e^{\sqrt{b_1-\lambda}})+\xi_4^{-^2}(e^{-\sqrt{b_1-\lambda}}-e^{-2\sqrt{b_1-\lambda}})\right],
\end{split}
\eeq
and
\beq
\begin{split}
I_2&= (a_1-a_2)\zeta_1^+\zeta_2^+ + \frac{a_1-a_2}{2\sqrt{a_2-\lambda}}\left[\zeta_1^{+^2}(e^{\sqrt{a_2-\lambda}}-1)+\zeta_2^{+^2}(1-e^{-\sqrt{a_2-\lambda}})\right]\\
& +  (b_1-b_2)\zeta_3^+\zeta_4^+ + \frac{b_1-b_2}{2\sqrt{b_2-\lambda}}\left[\zeta_3^{+^2}(e^{2\sqrt{b_2-\lambda}}-e^{\sqrt{b_2-\lambda}})+\zeta_4^{+^2}(e^{-\sqrt{b_2-\lambda}}-e^{-2\sqrt{b_2-\lambda}})\right],
\end{split}
\eeq
where $\zeta^+$ is the same as $\xi^+$ in \eqref{E:xi_plus} with $(a,b,\kappa,s)$ replaced by $(a_2,b_2,\kappa_2,1/2)$, $\kappa_2$ is the same as $\kappa$ in \eqref{E:kappa} with $(a,b,s)$ replaced by $(a_2,b_2,1/2)$.

\begin{figure}[!ht]
  \begin{center}
  \scalebox{0.6}{\includegraphics{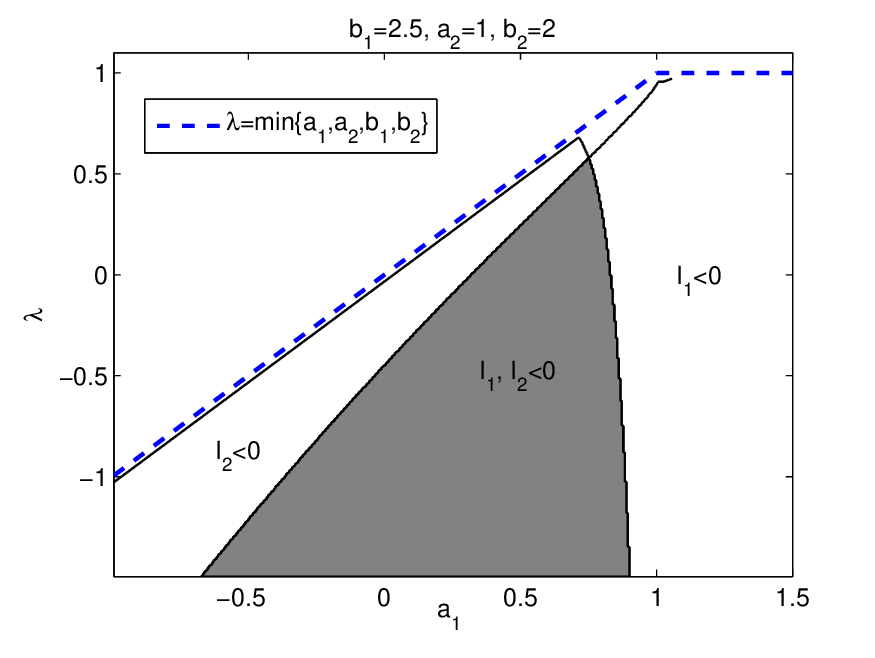}}
  \end{center}
  \caption{Region in the $(a_1,\lambda)$ plane, where suprema of the enclosures of $I_1, I_2$ are negative for the interface \eqref{E:V_interf} with piecewise constant $V_1$ and $V_2$ and $b_1=2.5, a_2=2, b_2=2, s_1=s_2=\tfrac{1}{2}$ (computations performed in interval arithmetic).}
  \label{F:distinct_Vs_2D}
\end{figure}

Fig. \ref{F:distinct_Vs_2D} shows regions of the $(a_1,\lambda)$ plane where the integrals $I_1,I_2$ are negative. The shaded region is where both $I_1$ and $I_2$ are negative, i.e. where ground state existence is guaranteed irrespectively of the coefficients $\Gamma_1, \Gamma_2$ and of $p$ (within (H1)-(H3)). Similarly to Fig. \ref{F:disloc_2D} we covered the region $\{(a_1,\lambda)\in \R^2 : -1\leq a_1\leq 1.5, -1.5 \leq \lambda \leq \min\{a_1,a_2,b_1,b_2\}-0.0256\}$ completely with squares of size $0.0056$ in each dimension and used interval arithmetic to evaluate $I_1$ and $I_2$.

\begin{figure}[!ht]
  \begin{center}
  \scalebox{0.55}{\includegraphics{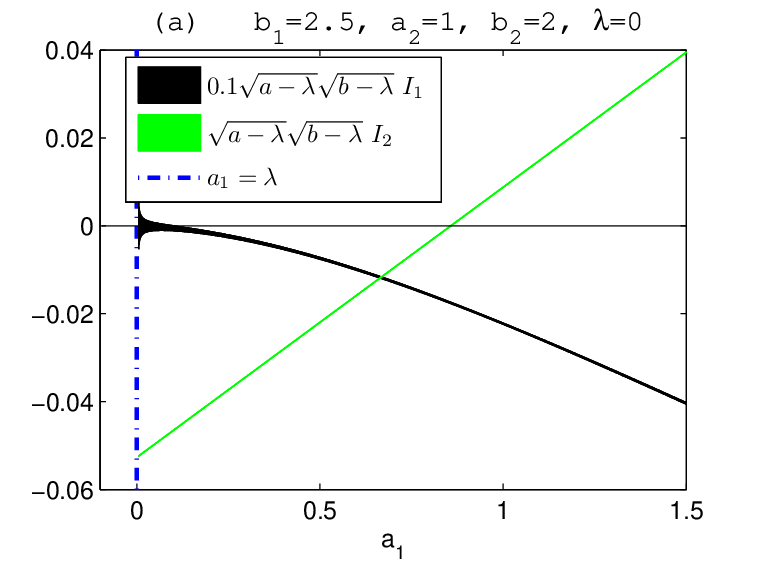}}
  \scalebox{0.55}{\includegraphics{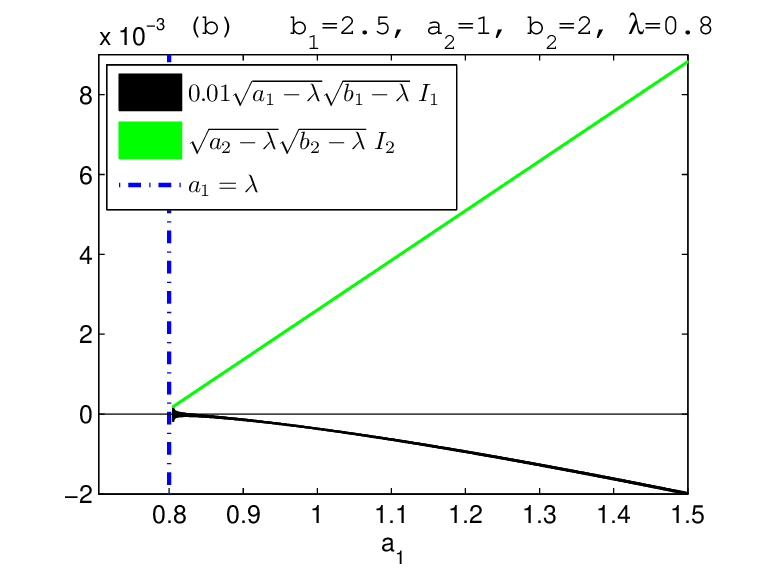}}
  \end{center}
  \caption{Scaled plots of the interval enclosures of $\sqrt{a_1-\lambda}\sqrt{b_1-\lambda} \ I_1$ and $\sqrt{a_2-\lambda}\sqrt{b_2-\lambda} \ I_2$ as functions of $a_1$ with $a_2,b_1$ and $b_2$ as in Fig. \ref{F:distinct_Vs_2D}. In (a) $\lambda=0$ and in (b) $\lambda = 0.8$ (computations performed in interval arithmetic with the $a_1$-interval width $1/1000$).}
  \label{F:distinct_Vs_b1_2p5}
\end{figure}

In Fig. \ref{F:distinct_Vs_b1_2p5} the enclosures of $\sqrt{a_1-\lambda}\sqrt{b_1-\lambda} \ I_1$ and $\sqrt{a_2-\lambda}\sqrt{b_2-\lambda} \ I_2$ as functions of $a_1$ are plotted for the example in Fig. \ref{F:distinct_Vs_2D} at two values of $\lambda$. At $\lambda=0.8$ the value of $I_2$ is always positive while at $\lambda=0$ both $I_1$ and $I_2$ are verified negative for $a_1\in (0.11, 0.85)$, where ground state existence thus follows.

\section{Piecewise Linear Potentials $V_1$ and $V_2$}\label{S:pwise_lin}

Here we consider continuous piecewise linear functions as potentials. Unlike in the case of piecewise constant potentials in Section \ref{S:pwise_const}, explicit formulas for the Bloch modes $u_{\pm}$ are now generally not available. We compute the Bloch waves via the numerical enclosure method presented in \cite{nag09}. All presented results are therefore verified in a strict mathematical sense.

\begin{example}\label{ex:1}
Let us consider the dislocation interface (\ref{E:disloc_interf}) 
with $\tau_1=\tau_2=:\tau$ for the following 
1-periodic potential
\[ V_0(x)=
\left\{
   \begin{array}{ll}
     4x-1, & x \in [0.25,0.5], \\
     -4x+3, & x \in [0.5,0.75], \\
     0, & x \in [0,0.25] \cup [0.75,1],
   \end{array}
   \right. \] 
cf. Fig. \ref{PL1q}.
\begin{center}
\begin{figure}[htbp]
\includegraphics[scale=0.6]{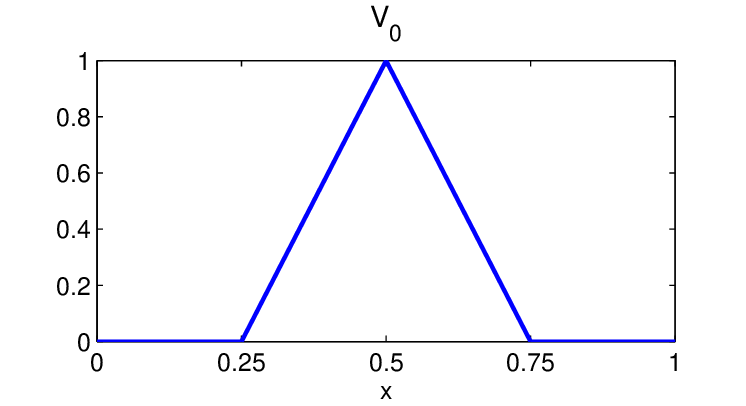}
\caption{Potential $V_0(x)$ on its period for example \ref{ex:1}}
\label{PL1q}
\end{figure}
\end{center}

Let $u_+$ and $u_-$ be a fundamental solution set of
\[ -u^{\prime \prime}+(V_0(x)-\lambda)u=0 \]
with the form
\begin{equation}
u_+(x)=e^{-\kappa x}p_{+}(x),~~u_-(x)=e^{\kappa x}p_{-}(x),
\end{equation}
where $\kappa>0$ is the characteristic exponent and $p_{+},p_{-}$
are $1-$periodic functions. 
Note that for $\rho:=e^{\kappa}>1$
\[u_+(x+1)=\rho^{-1} u_+(x),~~ u_-(x+1)=\rho u_-(x) \]
hold.

Since $V_0(x)$ is an even function, $u_+(-x)$ is also a solution.
$u_+(-x)$ must be a linear combination first of both $u_-(x)$ and 
$u_+(x)$, but since $u_+(x)$ grows at $-\infty$ whereas $u_+(-x)$ and 
$u_-(x)$ both decay at $-\infty$, the factor of $u_+(x)$ in the linear 
combination must be zero, i.e. $u_+(-x)=c u_-(x)$ holds for some
$c\in {\C}$.
Noting that the fundamental solutions can be normalized, 
we can define
\[ u_+(x):=u_-(-x) \]
after computing $u_-(x)$. 
Then by simple calculations we see that $I_1=I_2$ holds.

In Fig. \ref{PL1} the enclosure of $I_1$ as a function of $\tau$ is plotted for the two cases: (a) $\lambda=-1$, (b) $\lambda=-0.1$. In case of (a) $I_1$ is negative for $\tau \in [0.04,0.49]$ and in case of (b) $I_1$ is negative for $\tau \in [0.2,0.49]$. Hence ground states exist in these cases.

\begin{figure}[htbp]
\includegraphics[scale=0.6]{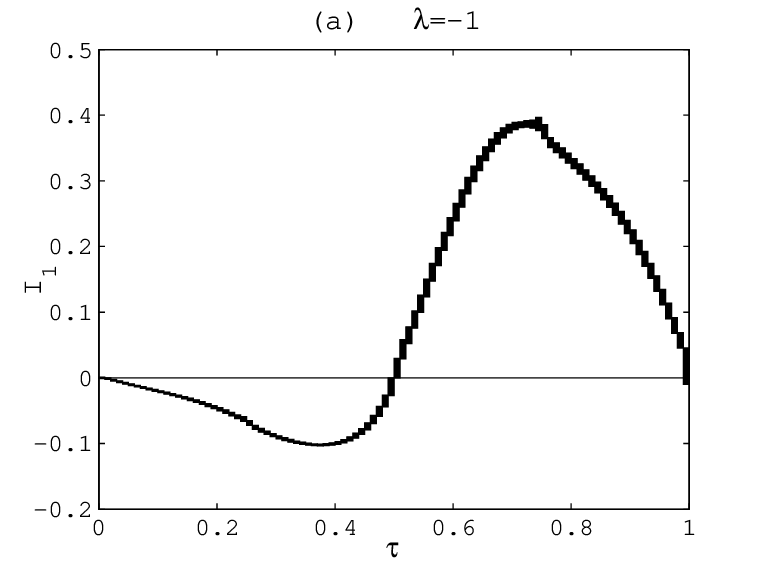}
\includegraphics[scale=0.6]{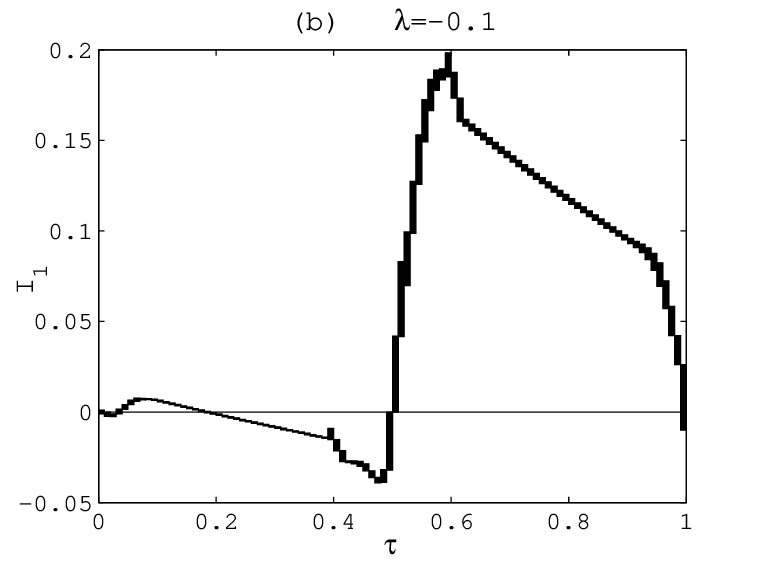}
\caption{Plot of the interval enclosure of $I_1$ as a function of $\tau$ with interval width 0.01 for example \ref{ex:1}.}
\label{PL1}
\end{figure}
\end{example}

\begin{example}\label{ex:2}
As another example of the dislocation interface (\ref{E:disloc_interf}) 
with $\tau_1=\tau_2=:\tau$ we consider the 1-periodic potential 
\[ V_0(x)=\left\{
   \begin{array}{ll}
     30x-4.5, & x \in [0.15,0.25], \\
     -30x+10.5, & x \in [0.25,0.35], \\
     0, & x \in [0,0.15] \cup [0.35,1]
   \end{array}
   \right. \] 

shown in Fig. \ref{PL2q}.
\begin{center}
\begin{figure}[htbp]
\includegraphics[scale=0.6]{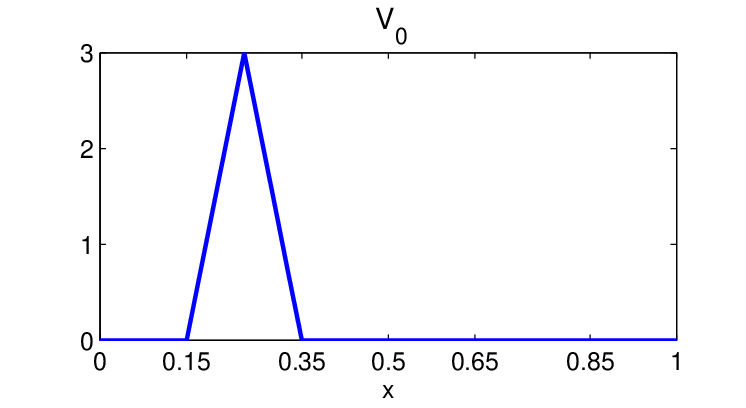}
\caption{Potential $V_0(x)$ on its period for example \ref{ex:2}}
\label{PL2q}
\end{figure}
\end{center}

Because $V_0$ is not even, $I_1$ is generally different from $I_2$ and we plot in Fig. \ref{PL2} the enclosure of both $I_1, I_2$ as functions of $\tau$ for the case (a) $\lambda=-1$ and (b) $\lambda=-0.01$. Several regions are observed, where $I_1$ or $I_2$ are negative and where by Theorem \ref{T:Bloch_cond} ground state existence follows for the corresponding interfaces  \eqref{E:disloc_interf}. As in Section \ref{S:pwise_const} we also find intervals where both $I_1$ and $I_2$ are negative so that ground states exist for arbitrary $\Gamma_1, \Gamma_2$ and $p$ (within H1-H3). In particular, in these intervals $\Gamma_1$ and $\Gamma_2$ do not need to be a dislocation of each other. 
For case (a) both $I_1$ and $I_2$ are negative for $\tau \in [0.01,0.31]\cup [0.51,0.66]$ and in case of (b) both $I_1$ and $I_2$ are negative for $\tau \in [0.18,0.31]$. 
\begin{figure}[htbp]
\includegraphics[scale=0.6]{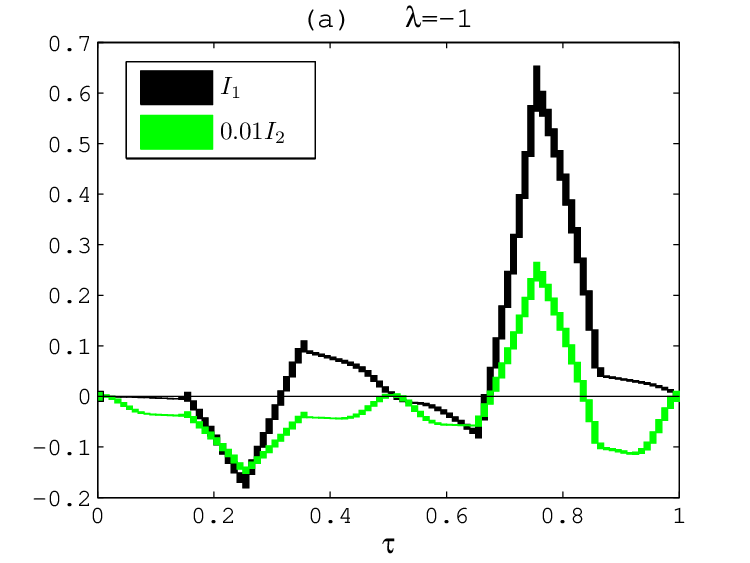}
\includegraphics[scale=0.6]{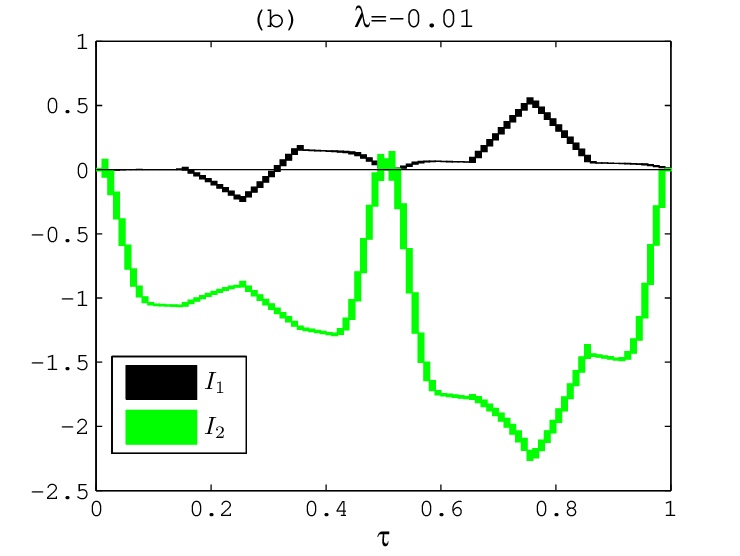}
\caption{Plot of the interval enclosure of $I_1, I_2$ as functions of $\tau$ with interval width 0.01 for example \ref{ex:2}.}
\label{PL2}
\end{figure}
\end{example}

\begin{example}\label{ex:3}
Here we consider a general interface (\ref{E:V_interf}) with piecewise linear $V_1, V_2$. As explained in Section \ref{S:pwise_const_gen_int}, we need to show that \textit{both} $I_1$ and $I_2$ are negative in order for Theorem \ref{T:Bloch_cond} to yield ground state existence.
As we noted below Corollary \ref{C:both_neg}, we should violate a monotone order between $V_1$ and $V_2$ in order to possibly obtain negative $I_1$ and $I_2$ simultaneously.

An example of potentials satisfying both $I_1<0$ and $I_2<0$ is
\begin{eqnarray}
V_1(x)&=&\left\{
   \begin{array}{ll}
     2000x-400, & x \in [0.2,0.25], \\
     -2000x+600, & x \in [0.25,0.3], \\
     0, & x \in [0,0.2] \cup [0.3,1],
   \end{array}\label{E:PL3q_V1}
   \right. \\ 
V_2(x)&=&\left\{
   \begin{array}{ll}
     20x-14, & x \in [0.7,0.75], \\
     -20x+16, & x \in [0.75,0.8], \\
     0, & x \in [0,0.7] \cup [0.8,1],
   \end{array}\label{E:PL3q_V2}
   \right.  
\end{eqnarray}
see Fig. \ref{PL3} (a).

\begin{center}
\begin{figure}[htbp]
\includegraphics[scale=0.7]{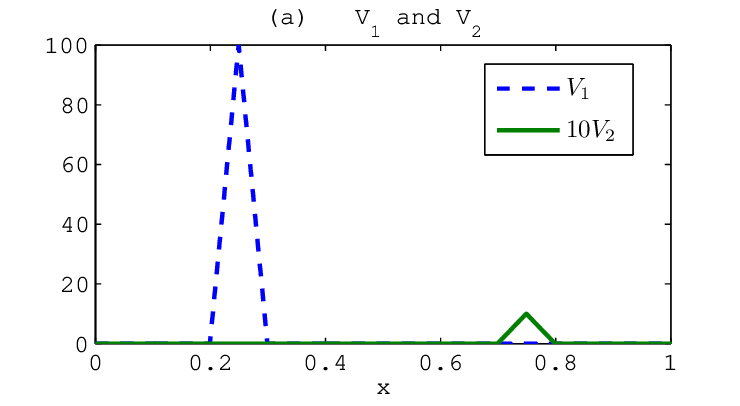}
\includegraphics[scale=0.6]{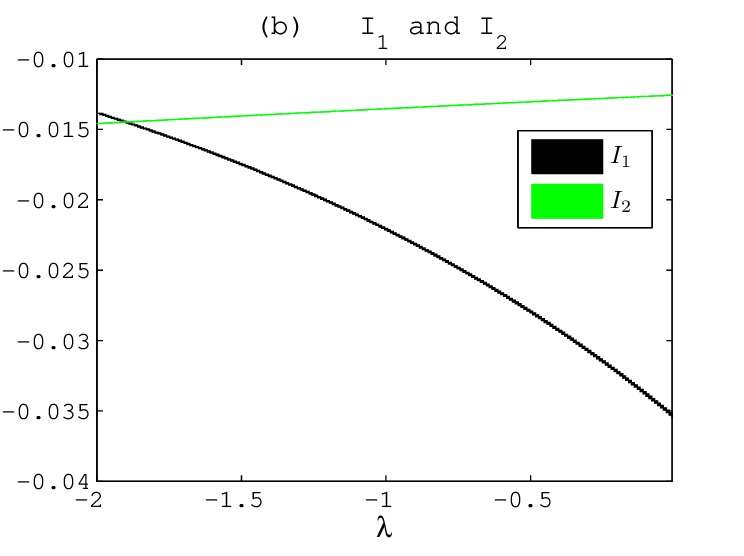}
\caption{(a) Potentials $V_1(x)$ and $V_2(x)$ in \eqref{E:PL3q_V1}, \eqref{E:PL3q_V2} on their period. (b) Interval enclosure of $I_1, I_2$ as functions of $\lambda$ with interval width 0.01 for potentials in  \eqref{E:PL3q_V1}, \eqref{E:PL3q_V2}.}
\label{PL3}
\end{figure}
\end{center} 

We have verified that for $\lambda\in [-2,-0.01]$ both $I_1$ and $I_2$ are negative. We plot in Fig. \ref{PL3} (b) the enclosure of both $I_1, I_2$ as functions of $\lambda$ using intervals with width 0.01.

In case when $\supp(V_2)$ is closer to $\supp(V_1)$ as given by \eqref{E:PL4q_V1}, \eqref{E:PL4q_V2} (see Fig. \ref{PL4q}), we verified that $I_1$ is 
negative and $I_2$ is positive for $\lambda\in [-2,-0.01]$, whence ground state existence cannot be concluded using Theorem \ref{T:Bloch_cond}.
\begin{eqnarray}
V_1(x)&=&\left\{
   \begin{array}{ll}
     2000x-400, & x \in [0.2,0.25], \\
     -2000x+600, & x \in [0.25,0.3], \\
     0, & x \in [0,0.2] \cup [0.3,1].
   \end{array}\label{E:PL4q_V1}
   \right. \\ 
V_2(x)&=&\left\{
   \begin{array}{ll}
     20x-6.5, & x \in [0.325,0.375], \\
     -20x+8.5, & x \in [0.375,0.425], \\
     0, & x \in [0,0.325] \cup [0.425,1].
   \end{array}\label{E:PL4q_V2}
   \right.  
\end{eqnarray}

\begin{center}
\begin{figure}[htbp]
\includegraphics[scale=0.7]{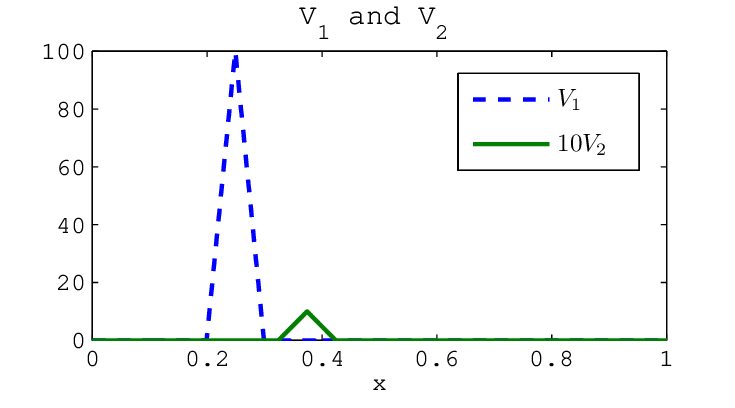}
\caption{Potentials $V_1(x)$ and $V_2(x)$ in \eqref{E:PL4q_V1} and \eqref{E:PL4q_V2} on their period.}
\label{PL4q}
\end{figure}
\end{center} 
\end{example}


\bibliographystyle{plain}
\bibliography{biblio_DNPR2011}


\end{document}